\documentstyle[12pt,multicol]{article}
\newtheorem{prop}{Proposition}[section]
\newtheorem{lemm}{Lemma}[section]
\newtheorem{coro}{Corollary}[section]
\newtheorem{defi}{Definition}[section]
\newtheorem{theo}{Theorem}[section]
\begin{document}

\title{Polyhedral Groups and Pencils of K3-Surfaces with Maximal Picard Number}
\author{W. Barth \and A. Sarti}
\date{\today}
\maketitle

\def \qed {\hspace*{\fill}\frame{\rule[0pt]{0pt}{8pt}\rule[0pt]{8pt}{0pt}}\par}
\def \qedup{\vskip-20pt\qed}

\def\R{{\rm I\!R}} 
\def\N{{\rm I\!N}} 
\def\F{{\rm I\!F}}
\def\P{{\rm I\!P}}
\def\C{{\mathchoice {\setbox0=\hbox{$\displaystyle\rm C$}\hbox{\hbox
to0pt{\kern0.4\wd0\vrule height0.9\ht0\hss}\box0}}
{\setbox0=\hbox{$\textstyle\rm C$}\hbox{\hbox
to0pt{\kern0.4\wd0\vrule height0.9\ht0\hss}\box0}}
{\setbox0=\hbox{$\scriptstyle\rm C$}\hbox{\hbox
to0pt{\kern0.4\wd0\vrule height0.9\ht0\hss}\box0}}
{\setbox0=\hbox{$\scriptscriptstyle\rm C$}\hbox{\hbox
to0pt{\kern0.4\wd0\vrule height0.9\ht0\hss}\box0}}}}

\def\Q{{\mathchoice {\setbox0=\hbox{$\displaystyle\rm Q$}\hbox{\raise
0.15\ht0\hbox to0pt{\kern0.4\wd0\vrule height0.8\ht0\hss}\box0}}
{\setbox0=\hbox{$\textstyle\rm Q$}\hbox{\raise
0.15\ht0\hbox to0pt{\kern0.4\wd0\vrule height0.8\ht0\hss}\box0}}
{\setbox0=\hbox{$\scriptstyle\rm Q$}\hbox{\raise
0.15\ht0\hbox to0pt{\kern0.4\wd0\vrule height0.7\ht0\hss}\box0}}
{\setbox0=\hbox{$\scriptscriptstyle\rm Q$}\hbox{\raise
0.15\ht0\hbox to0pt{\kern0.4\wd0\vrule height0.7\ht0\hss}\box0}}}}

\def\Z{{\mathchoice {\hbox{$\sf\textstyle Z\kern-0.4em Z$}}
{\hbox{$\sf\textstyle Z\kern-0.4em Z$}}
{\hbox{$\sf\scriptstyle Z\kern-0.3em Z$}}
{\hbox{$\sf\scriptscriptstyle Z\kern-0.2em Z$}}}}

\def\OO{{\cal O}}
\def\XX{\Xi}
\def\YY{{\Upsilon}}

\newcommand{\D}{\displaystyle} 
\newcommand{\vv}{\vspace{0.3cm}} 

\input prepictex
\input pictex
\input postpictex

\begin{abstract}

A K3-surface is a (smooth) surface which is simply connected and has trivial canonical bundle. In these notes we investigate three particular pencils of K3-surfaces with maximal Picard number. More precisely the general member in each pencil has Picard number $19$ and each pencil contains four surfaces with Picard number $20$. These surfaces are obtained as the minimal resolution of quotients $X/G$, where $G\subset SO(4,\R)$ is some finite subgroup and $X\subset \P_3(\C)$ denotes a $G$-invariant surface. The singularities of $X/G$ come from fix points of $G$ on $X$ or from singularities of $X$. In any case the singularities on $X/G$ are $A-D-E$ surface singularities. The rational curves which resolve them and some extra 2-divisible sets, resp. 3-divisible sets of rational curves generate the Neron-Severi group of the minimal resolution.

\end{abstract}

\begin{multicols}{2}
\tableofcontents
\end{multicols}

\setcounter{section}{-1}

\section{Introduction}

The aim of this note is to present three particular pencils of
K3-surfaces with Picard-number $\geq 19$. These three pencils are
related to the three polyhedral groups
$T,\, O,$ resp. $I$, 
(the rotation groups of the platonic solids tetrahedron, octahedron
and icosahedron) as follows: 
It is classical that the group $SO(4,\R)$ contains
central extensions
$$ \begin{array}{l|ccc}
      & G_6 & G_8 & G_{12}  \\ \hline
\mbox{of} & T \times T & O \times O & I \times I \\
 \end{array} $$
by $\pm 1$. Each group $G_n, \, n=6,8,12,$ has the obvious invariant
$q:=x_0^2+x_1^2+x_2^2+x_3^2$. In [S] it is shown that each group
$G_n$ admits a second non-trivial invariant
$s_n$ of degree $n$. (The existence of these invariants seems to have
been known before [Ra,C], but not their explicit form as computed
in [S].)
The pencil 
$$X_{\lambda} \subset \P_3(\C): \quad s_n+\lambda q^{n/2}=0$$
therefore consists of degree-$n$ surfaces admitting 
the symmetry group $G_n$. We consider here the pencil 
of quotient surfaces 
$$Y'_{\lambda}:= X_{\lambda}/G_n \subset \P_3/G_n.$$
It is - for us - quite unexpected that these (singular) surfaces 
have minimal resolutions $Y_{\lambda}$, 
which are K3-surfaces with Picard-number 
$\geq 19$.

In [S] it is
shown that the general surface $X_{\lambda}$ is smooth and that 
for each $n=6,8,12$ there are
precisely four singular surfaces $X_{\lambda}, \, \lambda \in \C$.
The singularities of these surfaces are ordinary nodes
(double points $A_1$) forming one orbit under $G_n$. 

For a smooth surface $X_{\lambda}$ the singularities on the quotient
surface $Y'_{\lambda}$ originate from fix-points of subgroups 
of $G_n$. Using [S, sect. 7] it is easy to enumerate these fix-points and to
determine the corresponding quotient singularities. On the minimal resolution 
$Y_{\lambda}$ of $Y'_{\lambda}$ we find enough rational curves
to generate a lattice in $NS(Y_{\lambda})$ of rank 19. In sect. 5 we show that
the minimal desingularisation $Y_{\lambda}$ is K3 and
that the structure of this surface varies with $\lambda$.
This implies that the general surface $Y_{\lambda}$ has Picard
number 19. Then in sect. 6.1 we use even sets [N], resp.
3-divisible sets [B, T] of rational curves to determine
completely the Picard-lattice of these surfaces $Y_{\lambda}$.

If $X_{\lambda}$ is one of the four nodal surfaces in the pencil,
there is an additional rational curve on $Y_{\lambda}$. This
surface then has Picard-number 20. (Such K3-surfaces usually
are called singular [SI].) We compute the Picard-lattice
for the surfaces $Y_{\lambda}$ in all twelve cases (sect. 6.2).

\section{Notations and conventions}

The base field always is $\C$. We abbreviate complex roots
of unity as follows:
$$\omega=e^{2 \pi i/3}=\frac{1}{2}(-1+\sqrt{-3}), \quad
\epsilon:=e^{2 \pi i/5}, \quad 
\gamma:=e^{2 \pi i/8}=\frac{1}{\sqrt{2}}(1+i).$$
By $G \subset SO(3)$ we always denote one of the (ternary) polyhedral groups
$T, \, O$ or $I$, and by $\tilde{G} \subset SU(2)$ the
corresponding binary group. By 
$$\sigma: SU(2) \times SU(2) \to SO(4)$$
we denote the classical $2:1$ covering. The group $G_n \subset SO(4),
n=6,8,12$, is the image $\sigma(\tilde{G} \times \tilde{G})$ for
$\tilde{G}=\tilde{T}, \tilde{O}, \tilde{I}$.
Usually we are interested more in the group
$$PG_n=G_n/\{\pm 1\} \subset PGL(4).$$
For $n=6,8,12$ it is isomorphic with $T \times T, O \times O, I \times I$
having the order $12^2=144, 24^2=576$, resp. $60^2=3600$.

\begin{defi} a) Let $id \not= g \in PG_n$. A {\em fix-line} for
$g$ is a line $L \subset \P_3$ with $gx=x$ for all $x \in L$.
The {\em fix-group} $F_L \subset PG_n$ is the subgroup consisting of all
$h \in PG_n$ with $hx=x$ for all $x \in L$. The {\em order} $o(L)$
of $L$ is the order of this group $F_L$.

b) The {\em stabilizer group} $H_L \subset PG_n$ is the subgroup 
consisting of all $h \in PG_n$ with $hL=L$. The {\em length}
$\ell(L)$ is the length
$$|PG_n|/|H_L|$$
of the $G_n$-orbit of $L$.

c) We shall encounter fix-lines of orders $2,3,4$ and $5$. We define
their {\em types} by
$$ \begin{array}{c|cccc}
\mbox{order} & 2 & 3 & 4 & 5 \\ \hline
\mbox{type}  & M & N & R & S \\
\end{array} $$
\end{defi}

\vv
We shall denote by $X_{\lambda}: s_n+\lambda q^{n/2}=0$ the
symmetric surface with parameter $\lambda \in \C$. All these surfaces
are smooth, but for four parameters $\lambda_i$.
These four singular parameters in the normalization of 
[S, p.445, p.449] are
$$ \begin{array}{cccc|cccc|cccc}
\multicolumn{4}{c|}{n=6} & \multicolumn{4}{c|}{n=8} &
                           \multicolumn{4}{c}{n=12} \\ \hline
\lambda_1 & \lambda_2 & \lambda_3 & \lambda_4 &
\lambda_1 & \lambda_2 & \lambda_3 & \lambda_4 &
\lambda_1 & \lambda_2 & \lambda_3 & \lambda_4  \\
-1 & -\frac{2}{3} & -\frac{7}{12} & -\frac{1}{4} &
-1 & -\frac{3}{4} & -\frac{9}{16} & -\frac{5}{9} &
-\frac{3}{32} & -\frac{22}{243} & -\frac{2}{25} & 0 \\
\end{array} $$
Sometimes we call the surface $X_{\lambda}$ of degree $n$ and
parameter $\lambda_i$ just $X_{n,i}$, or refer to it as
the case n,i.

\section{Fixpoints}

In this section we determine the fix-points for elements
$id \not=g \in PG_n$.

Recall 
that each $\pm 1 \not= p \in \tilde{G}$ has precisely two 
eigen-spaces in $\C^2$ with the product 
of its eigen-values $=det(p)=1$.

In coordinates $x_0,...,x_3$ on $\R^4$ the morphism
$\sigma:\tilde{G} \times \tilde{G} \to SO(4,\R)$
is defined by 
$\sigma(p_1,p_2):(x_k) \mapsto (y_k)$ with
$$ \left( \begin{array}{cc}
   y_0+iy_1 & y_2+iy_3 \\
  -y_2+iy_3 & y_0-iy_1 \\ \end{array} \right) =
p_1 \cdot \left( \begin{array}{cc}
               x_0+ix_1 & x_2+ix_3 \\
              -x_2+ix_3 & x_0-ix_1 \\ \end{array} \right) \cdot p_2^{-1}.$$
The quadratic invariant
$$q=x_0^2+x_1^2+x_2^2+x_3^2=det \left( \begin{array}{cc}
               x_0+ix_1 & x_2+ix_3 \\
              -x_2+ix_3 & x_0-ix_1 \\ \end{array} \right)$$
vanishes on tensor-product matrices
$$\left( \begin{array}{cc}
               x_0+ix_1 & x_2+ix_3 \\
              -x_2+ix_3 & x_0-ix_1 \\ \end{array} \right)
= \left( \begin{array}{cc}
          v_0w_0 & v_0w_1 \\ v_1w_0 & v_1w_1 \\
          \end{array} \right) = v \otimes w.$$
The action of $\tilde{G} \times \tilde{G}$ on the
quadric 
$$Q:=\{q=0\}=\P_1 \times \P_1$$ 
is induced by the
actions of the group $\tilde{G}$ on the tensor
factors $v$ and $w \in \C^2$  
$$\sigma(p_1,p_2):v \otimes w \mapsto (p_1v) \otimes (\bar{p_2}w).$$

The fix-points for $ \pm 1 \not= \sigma(p_1,p_2) \in G_n$ 
on $\P_3$ come in three kinds:

\vv
1) {\em Fix-points on the quadric:} $\pm 1 \not= p_1 \in \tilde{G}$ 
has two independent
eigenvectors $v,v'$. The spaces $v \otimes \C^2$ and $v' \otimes \C^2$
determine on the quadric two fix-lines 
for $\sigma(p_1,\pm 1)$ belonging to the same ruling. In this way 
$\tilde{G}$-orbits
of fix-points for elements $p_1 \in \tilde{G}$ determine $G_n$-orbits
of fix-lines in the same ruling of the following lengths: 
$$ \begin{array}{c|rrrr}
\mbox{order of }p & 4 & 6      & 8 & 10 \\ \hline
G_6    &            6 & 4,\, 4 & - & - \\
G_8    &           12 & 8      & 6 & - \\ 
G_{12} &           30 & 20     & - & 12 \\ \end{array} $$
In the same way fix-points for $p_2 \in \tilde{G}$ determine
fix-lines for $\sigma(\pm 1,p_2) \in G_n$ in the other ruling.
In [S, p.439] it is shown that the base locus of the pencil
$X_{\lambda}$ consists of $2n$ such fix-lines,
$n$ lines in each ruling, say $\Lambda_k,\Lambda_k',\, k=1,...,n$.
The fix-group $F_{\Lambda_k}$ for the general point on each line 
$\Lambda_k, \Lambda_k'$ then is cyclic
of order $s:=|G|/n$:
$$\begin{array}{c|rrr} n & 6 & 8 & 12 \\ \hline
                       s & 2 & 3 &  5 \\ \end{array} $$

Where a fix-line for $\sigma(p_1,\pm 1)$ meets a fix-line
for $\sigma(\pm 1, p_2)$ we obviously have an isolated fix-point $x$
for the group generated by these two symmetries. We denote by $t$
the order of the (cyclic) subgroup of $P(\sigma(\pm 1, \tilde{G}))$ fixing
$x$. The number of $H_{\Lambda_k}$-orbits on each line $\Lambda_k$
of such points is
$$ \begin{array}{rc|cccc}
   &   & \multicolumn{4}{c} t \\ 
n  & s & 2 & 3 & 4 & 5 \\  \hline
6  & 2 & 1 & 2 & - & - \\
8  & 3 & 1 & 1 & 1 & - \\
12 & 5 & 1 & 1 & - & 1 \\ \end{array} $$

\vv
2) {\em Fix-lines off the quadric:} Let $L \subset \P_3$ be a fix-line
for $\sigma(p_1,p_2) \in G_n$ with $p_1,p_2 \not= \pm 1$. It meets the quadric
in at least one fix-point defined by a tensor product
$v \otimes w$ with $v,w$ eigenvectors for $p_1,p_2$ respectively.
The group $<\sigma(p_1,\pm 1)> \subset H_L$ centralizes $\sigma(p_1,p_2)$.
Therefore there is a second fix-point on $L$ for this group.
Necessarily it lies on the quadric, being determined by a tensor-product
$v' \otimes w'$ with $v',w'$ eigenvectors for $p_1,p_2$ respectively.
Let $\alpha,\alpha'$ be the eigen-values for $p_1$ on $v,v'$ and
$\beta,\beta'$ those for $p_2$ on $w,w'$ respectively. Then
$$\alpha \cdot \alpha' = \beta \cdot \beta' = 1.$$
Since all points on $L$ have the same eigen-value under
$\sigma(p_1,p_2)$ we find
$$\alpha \cdot \beta = \alpha' \cdot \beta'=(\alpha \cdot \beta)^{-1}.$$
So $\alpha \cdot \beta= \pm 1$ and $g:=\sigma(p_1,p_2)$ acts on this line by 
an eigen-value $\pm 1$. 
In particular $p_1$ and $\pm p_2 \in \tilde{G}$
have the same order.

We reproduce from [S, p. 443] the table of $G_n$-orbits of fix-lines
off the quadric
by specifying a generator $g \in G_n$ of $F_L$. 
For this generator we use the notation of [S]. There it is also given
the length $\ell(L)$. This length determines the order 
$|H_L|=|PG_n|/\ell(L)$ of the
stabilizer group and the length $|H_L|/|F_L|$ of the general
$H_L$-orbit on $L$:
$$ \begin{array}{c|ccc|ccccc|ccc}
n           & \multicolumn{3}{c|}{6} &\multicolumn{5}{c|}{8} &
                                     \multicolumn{3}{c}{12} \\ \hline         
g           & \sigma_{24} & \pi_3\pi_3' & \pi_3^2 \pi_3' &
              \pi_3\pi_4\pi_3'\pi_4'& \pi_3\pi_4\sigma_4 &
              \sigma_2\pi_3'\pi_4' & \pi_3\pi_3' & \pi_4\pi_4' &
              \sigma_{24} & \pi_3\pi_3' & \pi_5\pi_5' \\
F_L         & \Z_2 & \Z_3 & \Z_3 & \Z_2 & \Z_2 & \Z_2 & \Z_3 & \Z_4 &
              \Z_2 & \Z_3 & \Z_5 \\
\mbox{type} & M  & N  & N' & M  & M' & M''& N  & R  & M   & N   & S  \\
\ell(L)     & 18 & 16 & 16 & 72 & 36 & 36 & 32 & 18 & 450 & 200 & 72 \\
|H_L|/|F_L| &  4 &  3 &  3 & 4  & 8  &  8 &  6 &  8 &  4  &  6  & 10 \\ 
\end{array} $$

\vv   
3) {\em Intersections of fix-lines off the quadric:}
From [S, p.450] one can read off the 
$G_n$-orbits of intersections of these
lines outside of the quadric and the value of the parameter
$\lambda$ for the surface $X_{\lambda}$ passing through
this intersection point. An intersection point is a
fix-point for the group generated by the transformations
leaving fixed the intersecting lines. In the following table
we give these (projective) groups ($D_n$ denoting the dihedral
group of order $2n$), the orders of the fix-group of intersecting
lines, the generators of these groups, 
as well as the numbers of lines meeting:
$$ \begin{array}{r|r|llll}
n & \lambda &
    \mbox{group} & \mbox{orders} & \mbox{generators} & \mbox{numbers} \\ \hline
6 & \lambda_1 &   T     &   2,3         & \sigma_{24}, \, \pi_3\pi_3'  &  3, 4    \\
  & \lambda_4 & T     &   2,3         & \sigma_{24}, \, \pi_3^2\pi_3'& 3, 4  \\  \hline
8 & \lambda_1  & O     &   2,3,4  &\pi_3\pi_4\pi_3'\pi_4',\,\pi_3\pi_3',\,\pi_4\pi_4'
                                                             & 6,4,3 \\
  & \lambda_2 & D_4   &   2,2,4  &\pi_3\pi_4\sigma_4,\, \sigma_2\pi_3'\pi_4', \,
                             \pi_4\pi_4'                     & 2,2,1 \\
  & \lambda_3 & \Z_2 \times \Z_2 & 2,2,2 &\pi_3\pi_4\sigma_4, \,\sigma_2\pi_3'\pi_4',
                            \, \pi_3\pi_4\pi_3'\pi_4'          & 1,1,1 \\
  & \lambda_4  &  D_3  & 2,3      & \pi_3\pi_4\pi_3'\pi_4',\,\pi_3\pi_3'& 3,1 \\
    \hline  
12 & \lambda_1 & T   & 2,3      & \sigma_{24}, \,\pi_3\pi_3'       & 3,4 \\
   & \lambda_2 &D_3 & 2,3      & \sigma_{24},\, \pi_3\pi_3'       & 3,1 \\
   & \lambda_3 & D_5 & 2,5      & \sigma_{24}, \,\pi_5\pi_5'       & 5,1 \\
   & \lambda_4 & I   & 2,3,5    & \sigma_{24}, \,\pi_3\pi_3',\, \pi_5\pi_5' 
                                                             & 15, 10, 6 \\
\end{array} $$

\section{Quotient singularities}

Singularities in the quotient surface $Y'=Y_{\lambda}'$ originate from
fix-points of the group action (or from singularities on $X$,
but the latter are included in the fix-points, see [S, (6.4)]).
We distinguish four types of fix-points on $X=X_{\lambda}$ for elements of
$G_n$:

\begin{itemize} 
\item[1)] Points of the base locus $\Lambda$ of the pencil, $n$
lines in each of the two rulings of the invariant quadric $Q$,
the (projective) fix-group being $\Z_s$ from section 1;
\item[2)] points on a line $\Lambda_k$ or $\Lambda_k'$ 
in the base locus, fixed by the
group $\Z_s=<\sigma(p,1)>$ from section 1 and by some
non-trivial subgroup $\Z_t \subset P(\sigma(1,\tilde{G}))$;
\item[3)] isolated fixed points on the intersection of a fix-line 
and a smooth surface $X_{\lambda}$; 
\item[4)] nodes of a surface $X_{\lambda}$.
\end{itemize}

{\bf 1)} All points of $\Lambda_i$
are fixed by the cyclic group $\Z_s$ from section 1.
The quotient map here is a cyclic covering of order $s$.
The quotient by $\Z_s$ is smooth.

\vv
{\bf 2)} Since $G_n$ acts on $\Lambda_i$ as the ternary polyhedral
group $G$, there are orbits of points on $\Lambda_i$,
fixed under some none-trivial subgroup of $G$. We have
to disinguish two cases:

{\em Case 1: The $n$ points, where the line $\Lambda_i$ meets some
line $\Lambda_k' \subset \Lambda$.} 
Here the stabilizer group is $\Z_s \times \Z_s$
acting on $X$ by reflections in the two lines 
$\Lambda_i,\Lambda_k'$.
In such points the quotient surface $Y'$ is smooth.

{\em Case 2: The fix-points of other non-trivial subgroups of $G$.}
The lengths of these orbits and their stabilizer subgroups 
$\Z_t \subset G$ are given in section 1:
$$ \begin{array}{c|rrr}
t       & 2    &      3 &  4    \\ \hline
G_6     &  -   & 4,\, 4 & -     \\
G_8     & 12   &   -    & 6     \\
G_{12}  & 30   &  20    & -     \\ 
\end{array} $$

The total stabilizer is the direct product $\Z_s \times \Z_t$.
Let $v,v'$ be eigen-vectors for $\Z_s$ and
$w,w'$ eigen-vectors for $\Z_t$. Let $v \otimes w$ determine
the fix-point in question. The surface $X$ is smooth there, containing the
line $\P(v \otimes \C^2)$, and intersecting the quadric 
$Q$ transversally. This implies that the tangent space of $X$ is
the plane
$$y_0 \cdot v \otimes w + y_1 \cdot v \otimes w' + y_2 \cdot v' \otimes w',
\quad y_0,y_1,y_2 \in \C.$$
Let $\sigma(p_1,\pm 1) \in \Z_s$ and $\sigma(\pm 1, p_2) \in \Z_t$ be
generators. Let them act by
$$\sigma(p_1,1)v=\alpha v,\, \sigma(p_1,1)v'=\alpha^{-1}v',\,
\sigma(1,p_2)w=\beta w,\, \sigma(1,p_2)w' = \beta^{-1} w'.$$
These transformations act on the coordinates $y_{\nu}$ of the
tangent plane as
$$ \begin{array}{l|ccccc}
            & y_0  & y_1 & y_2 & z_1:=y_1/y_0 & z_2:=y_2/y_0 \\ \hline
\sigma(p_1,1) & \alpha & \alpha & \alpha^{-1} & 1 & \alpha^{-2} \\
\sigma(1,p_2) & \beta  & \beta^{-1} & \beta^{-1} & \beta^{-2} &
                                                 \beta^{-2} \\
\end{array} $$

We introduce local coordinates on $X$ in which the group acts
as on $z_1,z_2$, and in fact use again $z_1,z_2$ to denote these
local coordinates on $X$. We locally form the quotient
$X/(\Z_s \times \Z_t)$ dividing first by the action of $\Z_s$
$$(z_1,z_2) \mapsto (z_1^s,z_2).$$
Then we trace the action of $\Z_t$ on $z_1^s$ and $z_2$.
A generator $\sigma(1,p_2)$ of $\Z_t$ acts by
$$ \begin{array}{ccccc}
y_0 & y_1 & y_2 & z_1 & z_2 \\ \hline
\omega & \omega^2 & \omega^2 &  \omega & \omega \\ 
 i     & -i       & -i       & -1      & -1     \\ 
\gamma & \gamma^7 & \gamma^7 & -i      & -i     \\ 
\end{array} $$
The resulting
singularities on $Y'$ are
$$\begin{array}{r|c|cccc|c}
n  & s & z_1 & z_2 & z_1^s & z_2 & \mbox{quotient singularity} \\ \hline
6  & 2 & \omega & \omega & \omega^2 & \omega & A_2 \\ \hline
8  & 3 & -1     & -1     &  -1      & -1     & A_1 \\
   &   & -i     & -i     &   i      & -i     & A_3 \\ \hline
12 & 5 & -1     & -1     & -1       & -1     & A_1 \\
   & 5 & \omega & \omega & \omega^2 & \omega & A_2 \\
\end{array} $$

{\bf 3)} Let $L \subset \P_3$ be a fix-line for $\sigma(p_1,p_2) \in G_n$,
not lying on the quadric. Assume that $\sigma(p_1,p_2)$ is chosen
as a generator for the group $F_L$. 
By sect. 1 there are eigen-vectors
$v,v'$ for $p_1$ with eigen-values $\alpha,\alpha^{-1}$ and
$w,w'$ for $p_2$ with eigen-values $\beta,\beta^{-1}$ 
satisfying 
$$\alpha \beta = \pm 1, \quad \alpha \beta = \alpha^{-1} \beta^{-1}
= \pm 1,$$ 
such that
$L$ is spanned by $v \otimes w$ and $v' \otimes w'$.
The general surface
$X$ meets this line in $n$ distinct points. If the line
has order $s$, two of these points lie on the base locus $\Lambda$.
So the number of points not in the quadric $Q$ cut out on $L$
by $X$ is
$$ \begin{array}{c|rr|rrr|rrr}
n             & \multicolumn{2}{c|}{6} & \multicolumn{3}{c|}{8} &
                                         \multicolumn{3}{c}{12} \\ \hline
o(L)          & 2 & 3 & 2 & 3 & 4 & 2 & 3 & 5 \\
\mbox{number} & 4 & 6 & 8 & 6 & 8 & 12 & 12 & 10 \\
\end{array} $$
These points fall into orbits under the stabilizer group $H_L$.
The lengths of these orbits are given in sect. 1.

To identify the quotient singularity we have to trace the action
of $\sigma(p_1,p_2)$ on the tangent plane $T_x(X)$. For general $X$
this plane will be transversal to $L$. So it must be the plane
spanned by $x, v \otimes w', v' \otimes w$. By continuity this
then is the case also for all smooth $X$. In particular, all smooth $X$ meet 
$L$ in $n$ distinct points, i.e., the intersections always are transversal.
And by continuity again, the numbers and lengths of $H_L$-orbits 
in $X \cap L$ are the
same for all smooth $X$. 
Since $\sigma(p_1,p_2)$ acts
$$ \begin{array}{c|ll}
\mbox{on} & v \otimes w' & v' \otimes w \\ \hline
\mbox{by} & \alpha \beta^{-1} & \alpha^{-1} \beta, \\
\end{array} $$
the eigen-values for $\sigma(p_1,p_2)$ on $T_x(X)$ are
$$ \frac{\alpha^{-1} \beta}{\alpha \beta} = \alpha^{-2} \quad
\mbox{ and } \quad \frac{\alpha \beta^{-1}}{\alpha \beta} =
\beta^{-2} = (\pm \frac{1}{\alpha})^2=\alpha^2. $$
The resulting quotient singularity on $Y'$ therefore is of
type $A_r$, where $r$ is the order of $L$.

We collect the results in the following table. It shows in
each case length and number of $H_L$-orbits, the number and type(s)
of the quotient singularity(ies). 
$$ \begin{array}{c|rrr|rrrrr|rrr}
n                    & \multicolumn{3}{c|}{6} & \multicolumn{5}{c|}{8}
                                              & \multicolumn{3}{c}{12} \\
                                                \hline
o(L)                 & 2 & 3 & 3 & 2 & 2 & 2 & 3 & 4 & 2 & 3 & 5 \\
\mbox{type}          & M & N'&N''& M'&M''& M & N & R & M & N & S \\
\mbox{length}        & 4 & 3 & 3 & 8 & 8 & 4 & 6 & 8 & 4 & 6 & 10 \\ 
\mbox{number}        & 1 & 2 & 2 & 1 & 1 & 2 & 1 & 1 & 3 & 2 & 1 \\
\mbox{singularities} & A_1 & 2A_2 & 2A_2 & A_1 & A_1 & 2A_1 & A_2 & A_3 &
                                          3A_1 & 2A_2 & A_4 \\ 
\end{array} $$

\vv
{\bf 4)} Finally we consider the nodal surfaces $X$. All the intersections of 
fix-lines considered in sect. 2 are nodes on the surfaces $X$. 
There are just two invariant surfaces
with nodes not given there, because through their nodes
passes just one fix-line. They are $G_6$-invariants. Their
parameters are as follows:
$$ \begin{array}{l|ll}
\lambda & \mbox{group} & \mbox{generator} \\ \hline
\lambda_2  & \Z_3        & \pi_3\pi_3' \\
\lambda_3  & \Z_3        & \pi_3\pi_3'^2 \\ \end{array} $$

We use this to collect the data for the twelve singular
surfaces $X$ in the next table. We include the number $ns$ of
nodes on the surface 
and specify the 
group $F \subset PSL(4)$ fixing the node.
For each type  
we give the number of lines meeting in the node. So
e.g. $3M$ means that there are three lines of type $M$ meeting at 
the node.
$$ \begin{array}{c|rrrr|rrrr|rrrr}
n       & 6  &      &       &      & 8  &      &       &      &
                                     12 &      &       &       \\ \hline
\lambda & \lambda_1 & \lambda_2 & \lambda_3  & \lambda_4 &
          \lambda_1 & \lambda_2 & \lambda_3  & \lambda_4 &
          \lambda_1 & \lambda_2 & \lambda_3  & \lambda_4 \\
ns      & 12 & 48   & 48    & 12   & 24 & 72   & 144              & 96  &
                                    300 & 600  & 360              & 60  \\    
F       & T  & \Z_3 & \Z_3  & T    & O  & D_4  & \Z_2 \times \Z_2 & D_3 &
                                     T  & D_3  & D_5   & I \\ \hline
        & 3M & 1N' & 1N'' & 3M &
          6M & 2M' & 1M' & 3M &
          3M & 3M & 5M & 15M \\
        & 4N' &    &     & 4N'' &  
          4N  & 2M''& 1M'' & 1N &
          4N  & 1N &  1S & 10N \\
        &        &    &          &            &
          3R    & 1R & 1M & & & & & 6S \\
\end{array} $$

\begin{lemm} Let $G \subset SO(3)$ be a finite subgroup of order $\geq 3$.

a) Up to $G$-equivariant linear coordinate change, there is a unique
$G$-invariant quadratic polynomial defining a non-degenerate
cone with top at the origin.

b) If $X$ is a $G$-invariant surface, 
having a node at the origin, then there is a
$G$-equivariant change of local (analytic) coordinates, such that
$X$ is given in the new coordinates by $x^2+y^2+z^2=0$. \end{lemm}

Proof. a) We distinguish two cases:

i) $G=\Z_2 \times \Z_2$ generated by the symmetries
 $$ (x,y,z) \mapsto (x,-y,-z) \mbox{ and }
(x,y,z) \mapsto (-x,y,-z).$$
The quadratic $G$-invariants are generated by the squares $x^2,y^2$ and
$z^2$. The invariant polynomial then is of the form
$ax^2+by^2+cz^2$ with $a,b,c \not= 0$. The coordinate change
$$x':=\sqrt{a}x, \, y':=\sqrt{b}y, \, z'=\sqrt{c}z$$
is $G$-equivariant and transforms the polynomial into
$x'^2+y'^2+z'^2$.

ii) $G$ contains an element $g$ of order $\geq 3$. Let it act by
$$(x,y,z) \mapsto (cx-sy, \, sx+cy, z)$$
with $c=cos(\alpha), \, s=sin(\alpha)$ and $\alpha \not= 0, \pi$.
The quadratic invariants of $g$ are generated by $x^2+y^2$ and 
$z^2$. The invariant polynomial must be of the form
$a(x^2+y^2)+bz^2$ with $a,b \not= 0$. The $G$-equivariant transformation
$x':=\sqrt{a}x,\, y':=\sqrt{a}y, \, z':=\sqrt{b}z$ transforms it 
into the same normal form as in i). This proves the assertion if
$G=<g>$ is cyclic or if $G$ is dihedral.

In the three other cases $G=T,O$ or $I$, it is well-known that
$x^2+y^2+z^2$, up to a constant factor, is the unique quadratic
$G$-invariant.

b) Let $X$ be given locally at the origin by an equation $f(x,y,z)=0$
with $f$ some power series. Since $X$ is $G$-invariant, so is the tangent cone
of $X$ at the origin. By a) we therefore may assume
$$f=x^2+y^2+z^2+f_3(x,y,z)$$
with a power series $f_3$ containing monomials of degrees $\geq 3$
only. It is well-known that there is a local biholomorphic
map $\varphi:(x,y,z) \mapsto (x',y',z')$ mapping $X$ to its tangent cone,
i.e., with the property $\varphi^*(x'^2+y'^2+z'^2)=f$. 
For the derivative $\varphi'(0)$ this implies
$\varphi'(0)^*(x'^2+y'^2+z'^2)=x^2+y^2+z^2$. 
After replacing $\varphi$ by $\varphi'(0)^{-1} \circ \varphi$
we even may assume $\varphi'(0)=id$.

Now consider the local
$G$-equivariant holomorphic map
$$\Phi:(x,y,z) \mapsto \frac{1}{|G|} \sum_{h \in G}
h \circ \varphi \circ h^{-1}.$$
Using the $G$-invariance of $f$ and $x'^2+y'^2+z'^2$ one easily checks
$\Phi^*(x'^2+y'^2+z'^2)=f$. It remains to show, that $\Phi$ locally
at the origin is biholomorphic. But this follows from
$$\Phi'(0)=\frac{1}{|G|} \sum_{h \in G} 
h \circ \varphi'(0) \circ h^{-1} =id.$$ \qed

\vv
Now consider the automorphism
$$\C^2 \to \C^2, \quad v=(v_0,v_1) \mapsto v^{\perp}:=(v_1,-v_0).$$
For $q \in SU(2)$ it is easy to check that $(qv)^{\perp} =\bar{q} v^{\perp}$.
Map $\C^2 \to \C^3$ via $v \mapsto v \otimes v^{\perp}$. Consider
$\C^3$ as the space of traceless complex matrices
$$\left( \begin{array}{cc}
          ix  &  y+iz \\ -y+iz & -ix \\ \end{array} \right).$$
Then 
$$v \otimes v^{\perp} = \left( \begin{array}{cc}
                         v_0v_1 & -v_0^2 \\ v_1^2 & -v_0v_1 \\
                        \end{array} \right)$$
is a matrix of determinant $x^2+y^2+z^2=0$. One easily checks that
the map $v \mapsto v \otimes v^{\perp}$ is $2:1$ onto the
cone of equation $x^2+y^2+z^2=0$, identifying this cone with
the quotient $\C^2/<-id>$. And this map is $SU(2)$-equivariant
with respect to the $2:1$ cover $SU(2) \to SO(3)$. 
If $\tilde{G} \subset SU(2)$ is some finite group, then
the quotient $\C^2/\tilde{G}$ via this map is identified with
the quotient of the cone by the corresponding ternary group 
$G \subset SO(3)$.

Together with
lemma 3.1 this shows:

\begin{prop} Let $X=X_{\lambda}$ be a nodal surface with $G$
the fix-group $G$ of the node. Then the image of
this node on $X/G_n$ is a quotient singularity locally isomorphic
with $\C^2/\tilde{G}$. \end{prop}

\section{Rational curves}

We denote by
$X=X_{\lambda} \to Y'=Y'_{\lambda}$
the quotient map for $G_n$ acting on $X$ and by
$Y=Y_{\lambda} \to Y'$
the minimal resolution of the quotient singularities on $Y$
coming from the orbits of isolated fixed points in sect. 2. 
The $n$ lines 
$\Lambda_i,\Lambda_i' \subset Q$
in each ruling map to one smooth rational curve in $Y'$.
We denote those by $L,L'$. Both 
these curves meet transversally in a smooth point of $Y'$.
All quotient singularities are rational double points.
Resolving them introduces more rational curves in $Y$.
For each singularity $A_k$ we get a chain of $k$ smooth rational (-2)-curves.
Since the group $\Z_t$ from sect. 3 acts on $X$ with $\Lambda_i$,
resp. $\Lambda_i'$ defining an eigen-space in the tangent space of $X$,
the curves $L,L'$ meet the $A_{t-1}$-string in an end curve of this
string, avoiding the other curves of the string.

All lines $L$ of the types $M,M',M'',N,N',N'',R,S$ form one orbit
under $G_n$. We denote by $M_i$ etc. the rational curves resolving
the $A_r$-singularity on the image of $L \cap X$. If
$L \cap X$ consists of more than one $H_L$-orbit we get in
this way more than one $A_r$-configuration of rational
curves coming from $L \cap X$. 
 
\subsection{The general case}

First we consider the quotients of the smooth surfaces $X$:
The striking
fact is that 
the number of the additional rational curves is 17. We give the dual 
graphs of the collections of 19 rational curves on $Y$, changing
the notation $L,L'$ to $L_3,L_3'$ for $n=6,12$ and to
$L_4,L_4'$ for $n=8$:
$$ 
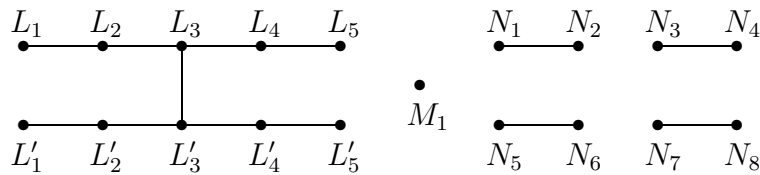
\beginpicture
\put {\line(1,0){120}} [Bl] at 0 30
\put {\line(1,0){120}} [Bl] at 0 0
\put {\line(1,0){30}} [Bl] at 180 30
\put {\line(1,0){30}} [Bl] at 180 0
\put {\line(1,0){30}} [Bl] at 240 30
\put {\line(1,0){30}} [Bl] at 240 0
\put {\line(0,1){30}} [Bl] at 60 0
\put {\circle*{4}} [Bl] at 0 30
\put {\circle*{4}} [Bl] at 30 30
\put {\circle*{4}} [Bl] at 60 30
\put {\circle*{4}} [Bl] at 90 30
\put {\circle*{4}} [Bl] at 120 30
\put {\circle*{4}} [Bl] at 0 0
\put {\circle*{4}} [Bl] at 30 0
\put {\circle*{4}} [Bl] at 60 0
\put {\circle*{4}} [Bl] at 90 0
\put {\circle*{4}} [Bl] at 120 0
\put {\circle*{4}} [Bl] at 150 15
\put {\circle*{4}} [Bl] at 180 30
\put {\circle*{4}} [Bl] at 210 30
\put {\circle*{4}} [Bl] at 240 30
\put {\circle*{4}} [Bl] at 270 30
\put {\circle*{4}} [Bl] at 180 0
\put {\circle*{4}} [Bl] at 210 0
\put {\circle*{4}} [Bl] at 240 0
\put {\circle*{4}} [Bl] at 270 0
\put {$L_1$} [Bl] at -5 35
\put {$L_2$} [Bl] at 25 35
\put {$L_3$} [Bl] at 55 35
\put {$L_4$} [Bl] at 85 35
\put {$L_5$} [Bl] at 115 35
\put {$L'_1$} [Bl] at -5 -15
\put {$L'_2$} [Bl] at 25 -15
\put {$L'_3$} [Bl] at 55 -15
\put {$L'_4$} [Bl] at 85 -15
\put {$L'_5$} [Bl] at 115 -15
\put {$M_1$} [Bl] at 145 0
\put {$N_1$} [Bl] at 175 35
\put {$N_2$} [Bl] at 205 35
\put {$N_3$} [Bl] at 235 35
\put {$N_4$} [Bl] at 265 35
\put {$N_5$} [Bl] at 175 -15
\put {$N_6$} [Bl] at 205 -15
\put {$N_7$} [Bl] at 235 -15
\put {$N_8$} [Bl] at 265 -15
\endpicture
$$
\begin{center}
n=6: $M_1$ coming from $M$, $N_1,...,N_4$ from $N'$, $N_5,...,N_8$ from $N''$
\end{center}

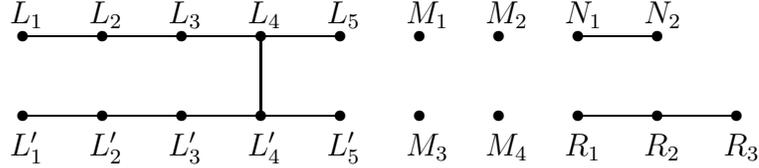
$$ 
\beginpicture
\put {\line(1,0){120}} [Bl] at 0 30
\put {\line(1,0){120}} [Bl] at 0 0
\put {\line(0,1){30}} [Bl] at 90 0
\put {\line(1,0){30}} [Bl] at 210 30
\put {\line(1,0){60}} [Bl] at 210 0
\put {\circle*{4}} [Bl] at 0 30
\put {\circle*{4}} [Bl] at 30 30
\put {\circle*{4}} [Bl] at 60 30
\put {\circle*{4}} [Bl] at 90 30
\put {\circle*{4}} [Bl] at 120 30
\put {\circle*{4}} [Bl] at 0 0
\put {\circle*{4}} [Bl] at 30 0
\put {\circle*{4}} [Bl] at 60 0
\put {\circle*{4}} [Bl] at 90 0
\put {\circle*{4}} [Bl] at 120 0
\put {\circle*{4}} [Bl] at 150 30
\put {\circle*{4}} [Bl] at 180 30
\put {\circle*{4}} [Bl] at 150 0
\put {\circle*{4}} [Bl] at 180 0
\put {\circle*{4}} [Bl] at 210 30
\put {\circle*{4}} [Bl] at 240 30
\put {\circle*{4}} [Bl] at 210 0
\put {\circle*{4}} [Bl] at 240 0
\put {\circle*{4}} [Bl] at 270 0
\put {$L_1$} [Bl] at -5 35
\put {$L_2$} [Bl] at 25 35
\put {$L_3$} [Bl] at 55 35
\put {$L_4$} [Bl] at 85 35
\put {$L_5$} [Bl] at 115 35
\put {$L'_1$} [Bl] at -5 -15
\put {$L'_2$} [Bl] at 25 -15
\put {$L'_3$} [Bl] at 55 -15
\put {$L'_4$} [Bl] at 85 -15
\put {$L'_5$} [Bl] at 115 -15
\put {$M_1$} [Bl] at 145 35
\put {$M_2$} [Bl] at 175 35
\put {$M_3$} [Bl] at 145 -15
\put {$M_4$} [Bl] at 175 -15
\put {$N_1$} [Bl] at 205 35
\put {$N_2$} [Bl] at 235 35
\put {$R_1$} [Bl] at 205 -15
\put {$R_2$} [Bl] at 235 -15
\put {$R_3$} [Bl] at 265 -15
\endpicture
$$
\begin{center}
n=8: $M_1$ from $M'$, $M_2$ from $M''$, $M_3,M_4$ from $M$, 
$N_i$ from $N$, $R_i$ from $R$
\end{center}

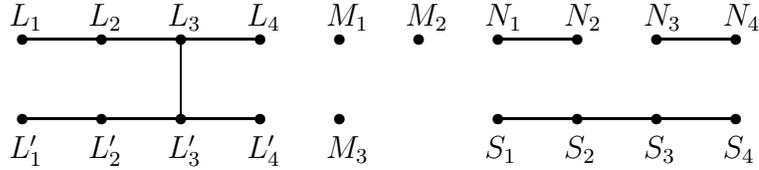
$$ 
\beginpicture
\put {\line(1,0){90}} [Bl] at 0 30
\put {\line(1,0){90}} [Bl] at 0 0
\put {\line(1,0){30}} [Bl] at 180 30
\put {\line(1,0){30}} [Bl] at 240 30
\put {\line(1,0){90}} [Bl] at 180 0
\put {\line(0,1){30}} [Bl] at 60 0
\put {\circle*{4}} [Bl] at 0 30
\put {\circle*{4}} [Bl] at 30 30
\put {\circle*{4}} [Bl] at 60 30
\put {\circle*{4}} [Bl] at 90 30
\put {\circle*{4}} [Bl] at 120 30
\put {\circle*{4}} [Bl] at 150 30
\put {\circle*{4}} [Bl] at 180 30
\put {\circle*{4}} [Bl] at 210 30
\put {\circle*{4}} [Bl] at 240 30
\put {\circle*{4}} [Bl] at 270 30
\put {\circle*{4}} [Bl] at 0 0
\put {\circle*{4}} [Bl] at 30 0
\put {\circle*{4}} [Bl] at 60 0
\put {\circle*{4}} [Bl] at 90 0
\put {\circle*{4}} [Bl] at 120 0
\put {\circle*{4}} [Bl] at 180 0
\put {\circle*{4}} [Bl] at 210 0
\put {\circle*{4}} [Bl] at 240 0
\put {\circle*{4}} [Bl] at 270 0
\put {$L_1$} [Bl] at -5 35
\put {$L_2$} [Bl] at 25 35
\put {$L_3$} [Bl] at 55 35
\put {$L_4$} [Bl] at 85 35
\put {$L'_1$} [Bl] at -5 -15
\put {$L'_2$} [Bl] at 25 -15
\put {$L'_3$} [Bl] at 55 -15
\put {$L'_4$} [Bl] at 85 -15
\put {$M_1$} [Bl] at 115 35
\put {$M_2$} [Bl] at 145 35
\put {$N_1$} [Bl] at 175 35
\put {$N_2$} [Bl] at 205 35
\put {$N_3$} [Bl] at 235 35
\put {$N_4$} [Bl] at 265 35
\put {$M_3$} [Bl] at 115 -15
\put {$S_1$} [Bl] at 175 -15
\put {$S_2$} [Bl] at 205 -15
\put {$S_3$} [Bl] at 235 -15
\put {$S_4$} [Bl] at 265 -15
\endpicture
$$

\begin{center}
n=12: $M_i$ from $M$, $N_i$ from $N$, $S_i$ from $S$
\end{center}

\begin{prop} In each case the 19 rational curves specified 
generate a sub-lattice
of $NS(Y)$ of rank 19. \end{prop} 

Proof. We compute
the discriminant $d$ of the lattice. The connected components of the dual graph
define sub-lattices, the direct sum of which is the lattice in question.
We compute the discriminant block-wise using the sub-lattices
$$ 
L := < L_i, L_i' >, \quad M:=<M_i>, \quad N:=<N_i>, \quad R:=<R_i>,
\quad S:=<S_i>$$
and find
$$ \begin{array}{r|rrrrr|l}
n  & d(L) & d(M) & d(N) & d(R) & d(S) & d \\ \hline
6  & -45   & -2  & 3^4  &      &      &  2 \cdot 3^6 \cdot 5 \\
8  & -28   & 2^4 & 3    & -4   &      & 2^8 \cdot 3 \cdot 7 \\
12 & -11   &-2^3 & 3^2  &      &  5   & 2^3 \cdot 3^2 \cdot 5 \cdot 11 \\
\end{array} $$ \qed

\subsection{The special cases}

Here we consider the desingularized quotients $Y$ for the twelve 
singular surfaces $X$. The image of the nodes on $X$ will be on $Y$
a quotient singularity for the binary group corresponding to the ternary group
$F$ from sect. 3. There we also gave the lines passing through this
node on $X$. The nodes of $X$ on such a line fall into orbits under the
group $H$ fixing the line. If there is just one $H$-orbit of intersection
points of the general surface $X$ with this line, it is clear that this orbit
converges to the orbit of nodes. We say: The quotient singularity
{\em swallows} the orbit. If however there are more than one
$H$-orbits, we have to analyze the situation more carefully.
We use the map onto $\P_1$ of this line induced by the parameter $\lambda$.
Nodes of $X$ on the given line will be branch points of this map.

\vv
{\em Degree 6}: On lines of type $M$ there ist just one orbit 
of four points. On lines of type $N',N''$ there are two orbits of length $3$.
The parameter $\lambda$ induces on each $N'$- or $N''$-line  
some $6:1$ cover over $\P_1$. Each fibre of six points
decomposes into two orbits of three points. The total ramification
degree is $-2 - 6\cdot (-2)=10$. The intersection with $Q$
consists of two points of ramification order $2$. So outside
of the quadric $Q$ we will have total ramification order six, hence
it will happen twice, that two orbits of three points are
swallowed by a quotient singularity. This must happen on
the surfaces $X_{6,1}$ and $X_{6,2}$ for $N'$, and for $N''$
on $X_{6,3}$ and $X_{6,4}$. We give the rational curves from 4.1 
disappearing in $Y$, being replaced by rational curves in the minimal 
resolution of the quotient surface. Here we do not mean that e.g. the curve
$N_1$ indeed converges to the curve denoted by $N_1$ in the dual graph
of the resolution. We just mean that all the curves denoted by
letters in the dual graph disappear:
$$
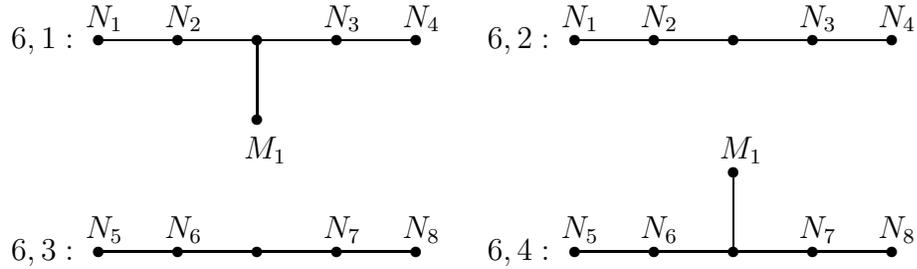
\beginpicture
\put {$6,1:$} [Bl] at -33 -3
\put {\line(1,0){120}} [Bl] at 0 0
\put {\line(0,1){30}} [Bl] at 60 -30
\put {\circle*{4}} [Bl] at 0 0
\put {\circle*{4}} [Bl] at 30 0
\put {\circle*{4}} [Bl] at 60 0
\put {\circle*{4}} [Bl] at 90 0
\put {\circle*{4}} [Bl] at 120 0
\put {\circle*{4}} [Bl] at 60 -30
\put {$N_1$} [Bl] at -5 5
\put {$N_2$} [Bl] at 25 5
\put {$N_3$} [Bl] at 85 5
\put {$N_4$} [Bl] at 115 5
\put {$M_1$} [Bl] at 55 -45
\put {$6,2:$} [Bl] at 147 -3
\put {\line(1,0){120}} [Bl] at 180 0
\put {\circle*{4}} [Bl] at 180 0
\put {\circle*{4}} [Bl] at 210 0
\put {\circle*{4}} [Bl] at 240 0
\put {\circle*{4}} [Bl] at 270 0
\put {\circle*{4}} [Bl] at 300 0
\put {$N_1$} [Bl] at 175 5
\put {$N_2$} [Bl] at 205 5
\put {$N_3$} [Bl] at 265 5
\put {$N_4$} [Bl] at 295 5
\put {$6,3:$} [Bl] at -33 -83
\put {\line(1,0){120}} [Bl] at 0 -80
\put {\circle*{4}} [Bl] at 0 -80
\put {\circle*{4}} [Bl] at 30 -80
\put {\circle*{4}} [Bl] at 60 -80
\put {\circle*{4}} [Bl] at 90 -80
\put {\circle*{4}} [Bl] at 120 -80
\put {$N_5$} [Bl] at -5 -75
\put {$N_6$} [Bl] at 25 -75
\put {$N_7$} [Bl] at 85 -75
\put {$N_8$} [Bl] at 115 -75
\put {$6,4:$} [Bl] at 147 -83 
\put {\line(1,0){120}} [Bl] at 180 -80
\put {\line(0,1){30}} [Bl] at 240 -80
\put {\circle*{4}} [Bl] at 180 -80
\put {\circle*{4}} [Bl] at 210 -80
\put {\circle*{4}} [Bl] at 240 -80
\put {\circle*{4}} [Bl] at 270 -80
\put {\circle*{4}} [Bl] at 300 -80
\put {\circle*{4}} [Bl] at 240 -50
\put {$N_5$} [Bl] at 175 -75
\put {$N_6$} [Bl] at 205 -75
\put {$N_7$} [Bl] at 265 -75
\put {$N_8$} [Bl] at 295 -75
\put {$M_1$} [Bl] at 235 -45
\endpicture $$

\vv
{\em Degree 8:} The only lines with two $H$-orbits are those of type
$M$. The map to $\P_1$ there has degree eight and total ramification order
$14$. The intersection with $Q$ counts for two points with ramification 
order three each. So there will be total ramification of order eight 
off the quadric. The surface $X_{8,1}$ has $24 \cdot 6/72=2$
nodes on such a line, it swallows at least one orbit. The surface
$X_{8,3}$ has $144/72=2$ nodes too and swallows at least one orbit too.
The surface $X_{8,4}$ has $96 \cdot 3/72=4$ nodes and swallows
at least two orbits. Since the total branching order adds up to
at least $2+2+4=8$, the bounds for the numbers of orders in
fact are exact numbers. The dual graphs for the resolution
of quotient singularities and the curves swallowed are as follows:
$$
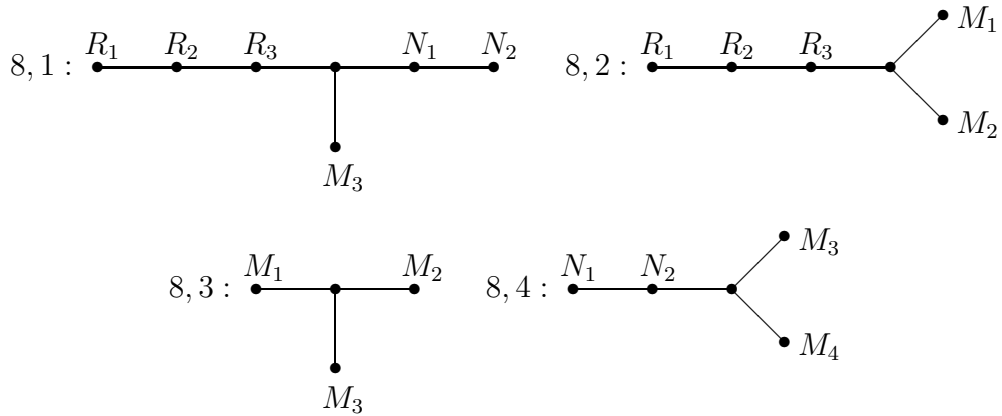
\beginpicture
\put {$8,1:$} [Bl] at -33 -3
\put {\line(1,0){150}} [Bl] at 0 0
\put {\line(0,1){30}} [Bl] at 90 -30
\put {\circle*{4}} [Bl] at 0 0
\put {\circle*{4}} [Bl] at 30 0
\put {\circle*{4}} [Bl] at 60 0
\put {\circle*{4}} [Bl] at 90 0
\put {\circle*{4}} [Bl] at 120 0
\put {\circle*{4}} [Bl] at 150 0
\put {\circle*{4}} [Bl] at 90 -30
\put {$R_1$} [Bl] at -5 5
\put {$R_2$} [Bl] at 25 5
\put {$R_3$} [Bl] at 55 5
\put {$N_1$} [Bl] at 115 5 
\put {$N_2$} [Bl] at 145 5 
\put {$M_3$} [Bl] at 85 -45
\put {$8,2:$} [Bl] at 177 -3
\put {\line(1,0){90}} [Bl] at 210 0
\put {\line(1,1){20}} [Bl] at 300 0
\put {\line(1,-1){20}} [Bl] at 300 0
\put {\circle*{4}} [Bl] at 210 0
\put {\circle*{4}} [Bl] at 240 0
\put {\circle*{4}} [Bl] at 270 0
\put {\circle*{4}} [Bl] at 300 0
\put {\circle*{4}} [Bl] at 320 20
\put {\circle*{4}} [Bl] at 320 -20
\put {$R_1$} [Bl] at 205 5
\put {$R_2$} [Bl] at 235 5
\put {$R_3$} [Bl] at 265 5
\put {$M_1$} [Bl] at 325 15
\put {$M_2$} [Bl] at 325 -25
\endpicture
$$

$$
\beginpicture
\put {$8,3:$} [Bl] at -33 -3
\put {\line(1,0){60}} [Bl] at 0 0
\put {\line(0,1){30}} [Bl] at 30 -30
\put {\circle*{4}} [Bl] at 0 0
\put {\circle*{4}} [Bl] at 30 0
\put {\circle*{4}} [Bl] at 60 0
\put {\circle*{4}} [Bl] at 30 -30
\put {$M_1$} [Bl] at -5 5
\put {$M_2$} [Bl] at 55 5
\put {$M_3$} [Bl] at 25 -45
\put{$8,4:$} [Bl] at 87 -3
\put {\line(1,0){60}} [Bl] at 120 0
\put {\line(1,1){20}} [Bl] at 180 0
\put {\line(1,-1){20}} [Bl] at 180 0
\put {\circle*{4}} [Bl] at 120 0
\put {\circle*{4}} [Bl] at 150 0
\put {\circle*{4}} [Bl] at 180 0
\put {\circle*{4}} [Bl] at 200 20
\put {\circle*{4}} [Bl] at 200 -20
\put {$N_1$} [Bl] at 115 5
\put {$N_2$} [Bl] at 145 5
\put {$M_3$} [Bl] at 205 15
\put {$M_4$} [Bl] at 205 -25
\endpicture
$$
Notice, that it is not necessary here to distinguish between
$M_3$ and $M_4$. In fact it is even impossible, since the two corresponding
orbits of intersections of the line $M$ with the surface 
$X_{\lambda}$ are interchanged by monodromy.

\vv
{\em Degree 12:} Now a line of type M contains three $H$-orbits of length four.
The total branching order for the $\lambda$-map is 22 on such a line.
The intersection with $Q$ consists of two six-fold points
and decreases the branching order by 10. So the total branching
order off the quadric is 12. On such a line there are
$$ \begin{array}{c|cc}
\mbox{on the surface }  & \mbox{nodes}        & \mbox{orbits swallowed}\\ \hline
 X_{12,1}               & 300 \cdot 3/450 = 2 & \geq 1 \\
 X_{12,2}               & 600 \cdot 3/450 = 4 & \geq 2 \\
X_{12,3}                & 360 \cdot 5/450 = 4 & \geq 2 \\
X_{12,4}                & 60 \cdot 15/ 450=2  & \geq 1 \\
\end{array}$$
Since the total branching order must add up to twelve, 
the number given is indeed the number of swallowed orbits.

A line of type $N$ contains two $H$-orbits of length six. Just as in the
preceding case one computes the following numbers
$$ \begin{array}{c|cc}
\mbox{on the surface} & \mbox{nodes} & \mbox{orbits swallowed}\\ \hline
X_{12,1}              & 300 \cdot 4/200 = 6 & \geq 2 \\
X_{12,2}              & 600 \cdot 1/200 = 3 & \geq 1 \\
X_{12,4}              & 60 \cdot 10/200 = 3 & \geq 1 \\
\end{array} $$
Again the total branching order adds up to twelve. Therefore the
estimates give the precise number of orbits swallowed.
$$ 
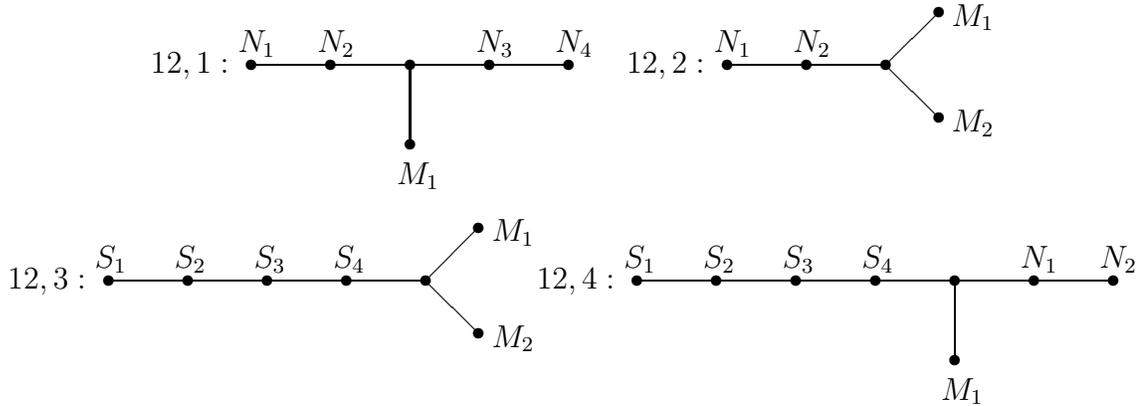
\beginpicture
\put {$12,1:$} [Bl] at -38 -3 
\put {\line(1,0){120}} [Bl] at 0 0
\put {\line(0,1){30}} [Bl] at 60 -30
\put {\circle*{4}} [Bl] at 0 0
\put {\circle*{4}} [Bl] at 30 0
\put {\circle*{4}} [Bl] at 60 0
\put {\circle*{4}} [Bl] at 90 0
\put {\circle*{4}} [Bl] at 120 0
\put {\circle*{4}} [Bl] at 60 -30
\put {$N_1$} [Bl] at -5 5
\put {$N_2$} [Bl] at 25 5
\put {$N_3$} [Bl] at 85 5
\put {$N_4$} [Bl] at 115 5
\put {$M_1$} [Bl] at 55 -45
\put {$12,2:$} [Bl] at 142 -3
\put {\line(1,0){60}} [Bl] at 180 0 
\put {\line(1,1){20}} [Bl] at 240 0
\put {\line(1,-1){20}} [Bl] at 240 0
\put {\circle*{4}} [Bl] at 180 0
\put {\circle*{4}} [Bl] at 210 0
\put {\circle*{4}} [Bl] at 240 0
\put {\circle*{4}} [Bl] at 260 20
\put {\circle*{4}} [Bl] at 260 -20
\put {$N_1$} [Bl] at 175 5
\put {$N_2$} [Bl] at 205 5
\put {$M_1$} [Bl] at 265 15
\put {$M_2$} [Bl] at 265 -25
\endpicture
$$

$$
\beginpicture
\put {$12,3:$} [Bl] at -38 -3 
\put {\line(1,0){120}} [Bl] at 0 0
\put {\line(1,1){20}} [Bl] at 120 0
\put {\line(1,-1){20}} [Bl] at 120 0
\put {\circle*{4}} [Bl] at 0 0
\put {\circle*{4}} [Bl] at 30 0
\put {\circle*{4}} [Bl] at 60 0
\put {\circle*{4}} [Bl] at 90 0
\put {\circle*{4}} [Bl] at 120 0
\put {\circle*{4}} [Bl] at 140 20
\put {\circle*{4}} [Bl] at 140 -20
\put {$S_1$} [Bl] at -5 5
\put {$S_2$} [Bl] at 25 5
\put {$S_3$} [Bl] at 55 5
\put {$S_4$} [Bl] at 85 5
\put {$M_1$} [Bl] at 145 15
\put {$M_2$} [Bl] at 145 -25
\put {$12,4:$} [Bl] at 162 -3
\put {\line(1,0){180}} [Bl] at 200 0
\put {\line(0,1){30}} [Bl] at 320 -30
\put {\circle*{4}} [Bl] at 200 0
\put {\circle*{4}} [Bl] at 230 0
\put {\circle*{4}} [Bl] at 260 0
\put {\circle*{4}} [Bl] at 290 0
\put {\circle*{4}} [Bl] at 320 0
\put {\circle*{4}} [Bl] at 350 0
\put {\circle*{4}} [Bl] at 380 0
\put {\circle*{4}} [Bl] at 320 -30
\put {$S_1$} [Bl] at 195 5
\put {$S_2$} [Bl] at 225 5
\put {$S_3$} [Bl] at 255 5
\put {$S_4$} [Bl] at 285 5
\put {$N_1$} [Bl] at 345 5
\put {$N_2$} [Bl] at 375 5
\put {$M_1$} [Bl] at 315 -45
\endpicture
$$
Again, by monodromy it is impossible to distinguish between
the curves $M_1,M_2$ and $M_3$, and likewise between the pairs
$\{N_1,N_2\}$ and $\{N_3,N_4\}$.

\section{K3-surfaces}

In this section we show that the 
desingularized quotient surfaces $Y_{\lambda}$ are K3 and that
their structure is not constant in $\lambda$.
We start with a crude but effective blow-up of $\P_3$. Let
$$\XX:=\{(x,\lambda) \in \P_3 \times \C: 
s_n(x)+\lambda q^{n/2}(x)=0\}.$$
In addition we put:

\begin{itemize}
\item $\bar{\XX} \subset \P_3 \times \P_1$ the closure of $\XX$.
It is a divisor of bidegree (n,1).
\item $\tau:\XX \to \P_3$ the natural projection onto the first factor;
\item $f:\XX \to \C$ the projection onto the second factor. It is given by
the function $\lambda$. 
\item $\tilde{\Lambda}:=\tau^{-1} \Lambda$. This pull-back of the base-locus
is the zero-set of $\tau^*q$ on $\XX$;
\item $\XX^0 \subset \XX$ the complement of the finitely many points
in $\XX$ lying over the nodes of the four nodal surfaces $X_{\lambda}$.  
\item $\YY' := \XX/G_n$ the quotient threefold. Notice that the
action of $G_n$ on $\P_3$ lifts naturally to an action on $\XX$.
\item $h:\YY' \to \C$ the map induced by $f$;
\item $\YY^0$ the image of $\XX^0$.
\end{itemize}

\begin{lemm} a) The threefold $\XX \subset \P_3 \times \C$ is smooth.

b) If $M \subset \P_3, \, M \not\subset Q,$ 
is a fix-line for an element $\pm 1 \not= g \in G_n$ 
and $\tilde{M} \subset \XX$ is its proper transform, then $\tilde{M}$
does not meet $\tilde{\Lambda}$ in $\XX$. 
\end{lemm}

Proof. a) By $\partial_{\lambda}(s_n+\lambda q^{n/2})= q^{n/2}$
singularities of $\XX$ can lie only on $\tau^{-1} \Lambda$. But
there
$$\partial_{x_i}(s_n+\lambda q^{n/2}) =\partial_{x_i} s_n.$$
Since $s_n=0$ is smooth along $\Lambda$, this proves 
that $\XX$ is smooth.

b) The assertion is obvious, if $M$ does not meet the base locus
$\Lambda$. If however $M \cap \Lambda = \{x_1,x_2\}$ is 
nonempty, we use the fact, observed in sect. 3, that 
the polynomial $s_n+tq^{n/2}|M$ vanishes in $x_i$ to the first
order for all
smooth surfaces $X: \, s_n+tq^{n/2}=0$. 
On $\tilde{M}$ however we have $s_n=-\lambda q^{n/2}$ with
$n/2 > 1$. So $\tilde{M}$ will not meet $ \tau^{-1}\{x_1,x_2\}$ 
in $\XX$. 
\qed

\vv
The $G_n$-action on $\XX$
has the following kinds of fix-points:

\begin{itemize}
\item[1)] Fix-points on $\tilde{\Lambda}$ for the group $\Z_s$;
\item[2)] Fix-points for the group $\Z_s \times \Z_s$ on the fibre
$\tau^{-1}(x)$ over some intersection of lines 
$\Lambda_i,\Lambda_j'$
in the base locus $\Lambda$;
\item[3)] Fix-points for a group $\Z_s \times \Z_t$ on the fibre
$\tau^{-1}(x)$ over a point $x$, where a line in the base locus
meets some line $M$ of fix-points not in the base locus. 
By lemma 5.1 b) $\tau^{-1}(x)$ and $\tilde{M}$ do not intersect
in $\XX$.
\item[4)] Fix-curves 
$\tilde{L}$ away from $\tilde{\Lambda}$ lying over fix-lines
$L$ not contained in the base-locus. All these curves are disjoint, when
considered in $\XX^0$.
\end{itemize}

The quotient three-fold $\YY'=\XX/G_n$ 
is smooth in the image points of fix-points of
types 1) or 2). It has quotient singularities $A_t$ in the image curves
of the curves $\tau^{-1}(x)$ of type 3). To be precise: The singularities there
locally are products of an $A_t$ surface singularity with a copy of the
complex unit disc.
Additional such cyclic quotient singularities $A_k$
occur on the image curves of curves $\tilde{L}$
of type 4). Where two such curves meet we have higher singularities.
But such points are removed in $\YY^0$. So $\YY^0$ is singular along finitely 
many smooth irreducible rational curves. The singularities along each curve
are products with some cyclic surface quotient $A_k$.

Let $\YY \to \YY^0$ be the minimal desingularisation of $\YY^0$
along these singular curves. Locally this is the product of the unit
disc with a minimal resolution of the surface singularity $A_k$.
Since the surfaces $Y_{\lambda}'$
intersect the singular curves transversally, the proper
transforms $Y_{\lambda} \subset \YY$ are smooth, minimally
desingularized. They are the fibres of the map induced by $h$.
For $\lambda_i, i=1,...,4,$ 
we denote by
$Y_{\lambda_i}$ the minimal resolutions of the quotient surfaces 
$X_{\lambda_i}/G_n$. We do not (and cannot) consider them
as surfaces in $\YY$.

\begin{prop} The surfaces $Y_{\lambda}$ are
(minimal) K3-surfaces. \end{prop}

Proof. All cyclic quotient singularities on $\YY^0$
are gorenstein. So there is a dualizing sheaf $\omega_{\YY^0}$
pulling back to the canonical bundle $K_{\YY}$ on $\YY$.
Under the quotient map $\XX^0 \to \YY^0$ it pulls back to the
canonical bundle $K_{\XX}$, except for points on the divisor
$\tilde{\Lambda}$. There we form the quotient in two steps, as
in sect. 3, first dividing by $\Z_s$ and then by $\Z_t$. The pull-back via
the quotient by $\Z_t$ is the canonical bundle of $\XX/\Z_s$. The
quotient map for $\Z_s$ is branched along $\tilde{\Lambda}$ to
the order $s$. So the adjunction formula shows: 
The dualizing sheaf $\omega_{\YY^0}$ pulls back to
$$K_{\XX^0} \otimes \OO_{\XX^0}((1-s) \tilde{\Lambda})
=K_{\XX^0} \otimes \tau^*(\OO_{\P_3}(2-2s)).$$
The divisor $\bar{\XX} \subset \P_3 \times \P_1$ is a divisor of 
bi-degree $(n,1)$. Hence $\bar{\XX}$ has a dualizing sheaf
$$\omega_{\bar{\XX}} = \OO_{\P_3 \times \P_1}(n-4,-2).$$
Now the miracle happens:
$$n-4=2s-2.$$
This implies: The pull-back of $\omega_{\YY^0}$ to $\XX^0$
equals the restriction of $\OO_{\P_3 \times \P_1}(0,-2)$,
i.e. it is trivial on $\XX^0$.

We distinguish two cases:

a) $\lambda \not= \lambda_i, \, i=1,...,4$: The adjunction formula for
$Y'=Y_{\lambda}'=X_{\lambda}/G_n$ shows
$$\omega_{Y'} = \omega_{\YY^0}|Y'.$$
So the pull-back of $\omega_{Y'}$ to $X$ is trivial.
This implies: $deg(\omega_{Y'})|C=0$
for all irreducible curves $C \subset Y'$ and then $deg(K_Y|C)=0$ for all
irreducible curves $C \subset Y$. The surfaces $Y$ have canonical
bundles, which are numerically trivial. In particular those surfaces are
all minimal. By the classification of algebraic surfaces
[BPV p. 168] they are abelian, K3, hyper-elliptic or Enriques.
Since we specified in sect. 6.1 rational curves on $Y$
spanning a lattice of rank 19 in $NS(Y)$ the only 
possibility is K3. 

b) $\lambda=\lambda_i, \, i=1,...,4$:
The proof of a) shows $deg(K_Y|C)=0$ for all
irreducible curves $C \subset Y$ not passing through the exceptional
locus of the minimal desingularization $Y \to Y'$. In particular
this holds for all curves $C$ which are proper transforms
of ample curves $D \subset Y'$. Now an arbitrary curve $C \subset Y$
is linearly equivalent to $E+C_1-C_2$ with $E$ exceptional
and $C_i$ proper transforms of ample curves $D_i \subset Y'$.
Since all singularities on $Y'$ are rational double points of
type A,D,E, we have $K_Y.E=0$.
The method from a) then applies here too. \qed

\begin{prop} The structure of the K3-surfaces $Y_{\lambda}$
varies with $\lambda$. \end{prop}

Proof. We restrict to surfaces near some surface $Y$, with
$Y'$ the quotient of a smooth surface $X$. Here we may 
assume that the total space $\YY$ is smooth. If all surfaces
near $Y$ were isomorphic, locally near $Y$ the fibration
would be trivial [FG]. I.e., there
would be an isomorphism $\Phi:Y \times D \to \YY$
respecting the fibre structure. Here $D$ is a copy
of the complex unit disc. By the continuity of the induced map
$$Y = Y \times{\lambda} \to Y_{\lambda}$$
there is an isomorphism $Y \to Y_{\lambda}$ mapping the 19
rational curves from sect. 4.1 on $Y$ to the corresponding curves on 
$Y_{\lambda}, \lambda \in D$. 
The covering $X \to Y'$ is defined by a subgroup in the fundamental
group of the complement in $Y$ of these rational curves. This
implies that the isomorphism $Y \to Y_{\lambda}$ induces
an isomorphism of the coverings $X \to X_{\lambda}$ equivariant
with respect to the $G_n$-action. 

Now this isomorphism must map the canonical bundle $\OO_X(n-4)$
to the canonical bundle $\OO_{X_{\lambda}}(n-4)$.
Since the surfaces $X_{\lambda}$ are simply-connected,
the isomorphism maps $\OO_X(1)$ to $\OO_{X_{\lambda}}(1)$,
i.e., it is given by a projectivity.  
This is in conflict with
the following. \qed

\begin{lemm} For general $\lambda \not= \mu$ there is no projectivity
$\varphi:\P_3 \to \P_3$ inducing some $G_n$-equivariant isomorphism
$X_{\lambda} \to X_{\mu}$.
\end{lemm}

Proof. Assume that such an isomorphism $\varphi$ exists.
Equivariance means for each $g \in G_n$ and $x \in X_{\lambda}$ that
$\varphi g (x)=g \varphi(x)$ or $\varphi^{-1}g^{-1} \varphi g (x) = x$. 
Since $X_{\lambda}$ spans $\P_3$ this implies the same property
for all $x \in \P_3$, i.e., the map $\varphi$ is $G_n$-equivariant
on all of $\P_3$. In particular, if $L \subset \P_3$ is a fixline
for $g \in G_n$, then so is $\varphi(L)$. Then we may as well assume
$\varphi(L)=L$.
We obtain a contradiction
by showing that the point sets $X_{\lambda} \cap L$ and
$X_{\mu} \cap L$ in general are not projectively equivalent.

{\em The cases n=6 and 12:}
We use the fix-line $L:=\{x_0=x_1=0\}$ of type $M$, fixed under
$\sigma_{1,3}= \sigma(q_1,q_1)$ (notation of [S, p. 432]). The group
$H_L$ has order 8, containing in addition the symmetries
$\sigma(q_1,1)$ and $\sigma(q_1q_2,q_1q_2)$ sending a point
$x=(0:0:x_2:x_3) \in L$ to 
$$\sigma(q_1,1)(x)=(0:0:x_2:-x_3), \quad
\sigma(q_1q_2,q_1q_2)(x)=(0:0:-x_3:x_2).$$
Omitting the first two coordinates and putting $x_2=1, x_3=u$, 
we find that a general $H_L$-orbit on $L$
consists of points
$$(1:u), \, (1:1/u), \, (1:-u), \, (1:-1/u).$$
The cross-ratio of these four points
$$CR = \frac{2u}{u+1/u}:\frac{1/u+u}{2/u}=
\frac{4u^2}{(1+u^2)^2}$$
varies with $u$. The intersection of $X_{6,\lambda}$ with $L$
consists of one such orbit, the intersection of $X_{12,\lambda}$
of three orbits. This implies the assertion for $n=6$ and 12.

{\em The case n=8:} Here we use the fix-line $L:=\{x_1=x_3,x_2=0\}$
of type $M$ for $\pi_3\pi_4\pi_3'\pi_4'$. Again $H_L$ has order 8
containing in addition $\pi_3\pi_4$ and $\sigma(q_1q_2,q_1q_2)$.
They send a point $x=(u:1:0:1) \in L$ to
$$\pi_3\pi_4(x)=(-2:u:0:u), \quad 
\sigma(q_1q_2,q_1q_2)(x)=(u:-1:0:-1).$$
Omitting the coordinates $x_3$ and $x_4$ we find that a general
$H_L$-orbit consists of
$$(u:1), \, (-u:1), \, (2/u:1), \, (-2/u:1).$$
Their cross-ratio
$$CR=\frac{u-2/u}{u+2/u}:\frac{-u-2/u}{-u+2/u}
=\frac{(u^2-2)^2}{(u^2+2)^2}$$
varies with $u$. The intersection of $X_{8,\lambda}$
consists of two such orbits. \qed

\begin{coro} The general K3-surface $Y_{\lambda}$ has
Picard-number 19. \end{coro}

\section{Picard-Lattices}

Here we compute the Picard lattices of our quotient K3-surfaces $Y$.

\subsection{The general case}

Denote by $V \subset H^2(X,\Z)$ the rank-19 lattice spanned 
(over $\Z$) by the 
rational curves from sect. 4.1.
For $n=6$ and $8$ this lattice $V$ 
is not the total Picard lattice:

\begin{prop} a) (n=6) After perhaps interchanging curves $N_{2i-1}$
and $N_{2i}$
the two divisor-classes
\begin{eqnarray*}
L  &:=& L_1-L_2+L_4-L_5+N_1-N_2+N_3-N_4+N_5-N_6+N_7-N_8, \\ 
L' &:=& L_1'-L_2'+L_4'-L_5'+N_1-N_2+N_3-N_4-N_5+N_6-N_7+N_8
\end{eqnarray*}
are divisible by 3 in NS(Y). Together with $V$
they span a rank-19 lattice with discriminant $2 \cdot 3^2 \cdot 5$.

b) (n=8) The two classes
\begin{eqnarray*}
L &:=& L_1+L_3+L_5+M_1+M_3+M_4+R_1+R_3, \\ 
L'&:=& L_1'+L_3'+L_5'+M_2+M_3+M_4+R_1+R_3
\end{eqnarray*}
are divisible by 2 in NS(Y). Together with $V$ 
they span a rank-19 lattice with discriminant
$2^4 \cdot 3 \cdot 7$. \end{prop}

Proof. a) Consider reduction modulo 3
$$\varphi_3: \Z^{22} =H^2(Y,\Z) \to H^2(Y,\F_3) = \F_3^{22}.$$
Because of
$$M_1^2=-2, \quad M_1.L_i'=0, \quad det(L_i',L_j')_{i,j=1,..,4}=5,$$
the images of $M_1,L_1',L_2',L_3',L_4'$ span a subspace of
$H^2(Y,\F_3)$ on which the intersection form has rank 5. The orthogonal
complement $C$ of this lattice in $H^2(Y,\F_3)$ has dimension 17
and the form is non-degenerate there. This $C$ 
contains the classes $mod \, 3$ of the twelve curves
$$L_1,L_2,L_4,L_5,N_1,...,N_8$$
Assume that
$$D_1:= \varphi_3<L_1,L_2,L_4,L_5,N_1,...,N_8>$$
has $\F_3$-dimension 12. Then
$$D_2:=\varphi_3<L_1-L_2,L_4-L_5,N_1-N_2,N_3-N_4,N_5-N_6,N_7-N_8>$$
has dimension six. Since $D_1 \perp D_2$, this is a contradiction.
We have shown: A non-trivial linear combination of the twelve classes
$L_1,L_2,L_4,L_5,N_1,...,N_8$ lies in the kernel of $\varphi_3$.
By [T] such a 3-divisible class contains at least 12 curves.
Hence we may assume the class is
$$L:=\lambda_1(L_1-L_2)+\lambda_4(L_4-L_5)+ 
\sum \nu_i(N_{2i-1}-N_{2i})$$
with $\lambda_j,\nu_i= \pm 1$ modulo 3. W.l.o.g. we put
$\lambda_1=1$. Intersecting with $L_3$ we find $\lambda_4=1$ too.
And after perhaps interchanging curves $N_{2i-1}$ 
with $N_{2i}$ we may assume
$\nu_1=...=\nu_4=1$. 

Exactly in the same way we find a class
$$L':=L_1'-L_2'+L_4'-L_5'+\sum \nu_i' (N_{2i-1}-N_{2i}), 
\quad \nu_i'=\pm 1 \, mod \, 3,$$
which is 3-divisible in $NS(Y)$. Then $L+L'$ is 3-divisible too,
and by [T] contains precisely 12 curves. This implies that
precisely two coefficients $\nu_i'$ cancel against the corresponding
coefficients of $L$. If these are the coefficients 
$\nu_3'$ and $\nu_4'$, we are done. If this should not be the case,
after perhaps interchanging $\{N_1,N_2\}$ with $\{N_3,N_4\}$,
$\{N_5,N_6\}$ with $\{N_7,N_8\}$ we may assume $\nu_1'=\nu_3'=1$
and $\nu_2'=\nu_4'=-1$. Denote by $T_2:H^2(X,\Z) \to H^2(X,\Z)$
the monodromy about $X_{6,2}$ (circling the parameter $\lambda_2$
in the parameter space) and by $T_3$ the monodromy about $X_{6,3}$.
So $T_2$ interchanges $\{N_1,N_2\}$ with $\{N_3,N_4\}$, leaving
fixed $\{N_5,N_6\},\{N_7,N_8\}$ with $T_3$ doing just the opposite.
$NS(X)$ contains the classes (coefficients modulo 3)
$$ \begin{array}{c|cccccc}
     & \frac{L_1+L_3+L_5}{3} & \frac{L_1'+L_3'+L_5'}{3} 
     & \frac{N_1-N_2}{3} & \frac{N_3-N_4}{3} 
     & \frac{N_5-N_6}{3} & \frac{N_7-N_8}{3} \\ \hline
L    &       1     &       0        &   1     &    1    &    1    &   1 \\
L'   &       0     &       1        &   1     &   -1    &    1    &  -1 \\
L+L' &       1     &       1        &  -1     &    0    &   -1    &   0 \\
T_2(L+L') &  1     &       1        &   0     & \pm 1   &   -1    &   0 \\
T_3(L+L') &  1     &       1        &  -1     &    0    &    0    & \pm 1 \\
\end{array} $$
These classes would span in $NS(X)/V$ a subgroup of order $3^4$,
in conflict with $d(V)=2 \cdot 3^6 \cdot 5$, contradiction.

b) Here we consider reduction modulo 2
$$\varphi_2:\Z^{22}=H^2(Y,\Z) \to H^2(Y ,\F_2) = \F_2^{22}.$$
The subspace
$$C:=\varphi_2 <L_1,L_3,L_5,M_1,M_2,M_3,M_4,R_1,R_3> \, \,
\subset H^2(Y,\F_2)$$
is totally isotropic. It is orthogonal to 
$D:=\varphi_2<L_1',L_2',L_3',L_4',N_1,N_2>$.
Because of
$$det (L_i'.L_j')_{i,j=1,...,4}=5, \quad det(N_i.N_j)_{i,j=1,2}=3,$$
the intersection form on $D$ is non-degenerate, and $D^{\perp}$
is non-degenerate of rank 16.
This implies $dim \, C \leq 8$. So there is a class
$$L:=\sum \lambda_i L_i+\mu_iM_i+\rho_iR_i$$
in the kernel of $\varphi_2$. By [N] it has precisely eight coefficients $=1$.
Intersecting 
$$ \begin{array}{l|l}
\mbox{with} & \mbox{we find} \\ \hline
L_2, L_4    & \lambda_1=\lambda_3=\lambda_5=:\lambda \\
R_2         & \rho_1=\rho_3=:\rho \\ \end{array} $$
This implies that precisely one coefficient $\mu_i$ will vanish and
$\lambda=\rho=1$. In the same way one finds a class 
$$L':=L_1'+L_3'+L_5'+\sum \mu_i'M_i+R_1+R_3$$
in the kernel of $\varphi_2$ with precisely one $\mu_1'$
vanishing. The class
$$L+L'=L_1+L_3+L_5+L_1'+L_3'+L_5'+\sum (\mu_i+\mu_i')M_i$$
also is divisible by $2$ and has precisely eight non-zero
coefficients. It follows that precisely two of the non-zero coefficients
from $\mu_i$ and $\mu_i'$ coincide.
If $\mu_3=\mu_4=\mu_3'=\mu_4'=1$ we are done (perhaps after
interchanging $L$ and $L'$). If this is not the case, assume e.g.
$\mu_1=\mu_2=\mu_4=1, \, \mu_3=0$. Denote by $T$ the monodromy
about the surface $X_{8,4}$ (circling the parameter $\lambda_4$
in the parameter space). It interchanges $M_3$ and $M_4$.
So there are the three classes
$$ \begin{array}{c|cccccc}
     & \frac{L_1+L_3+L_5}{2} & \frac{L_1'+L_3'+L_5'}{2} 
     & \frac{M_1}{2} &\frac{M_2}{2} & \frac{M_3}{2} & \frac{M_4}{2} \\
     \hline
L/2    &   1   &   0   &   1  &  1  &  0  &    1  \\
T(L/2) &   1   &   0   &   1  &  1  &  1  &    0  \\
L'/2   &   0   &   1   &   *  &  *  &  *  &    *  \\
\end{array} $$
belonging to $NS(X)$. Together with $V$ they span a lattice $W$
with discriminant $2^2 \cdot 3 \cdot 7$ and 
$[NS(X)^{\vee}:W] \leq 2^2 \cdot 3 \cdot 5$. However there are the two
classes
$$M:=M_1+M_2, \quad 
R:=R_1+2R_2+3R_3.$$
The two classes $M/2$ and $R/4$ belong to $NS(X)^{\vee}$
spanning in $NS(X)^{\vee}/W$ a subgroup $\simeq \Z_2 \times \Z_4$
of order $2^3$, contradiction.
\qed

\begin{theo} If the Neron-Severi group of $Y$ has rank 19, it
is generated by $V$ and
$$ \begin{array}{r|l}
n &           \\ \hline
6 & L/3, L'/3, \\
8 & L/2, L'/2, \\
12 & \mbox{no other classes.} \\
\end{array} $$
\end{theo}

Proof a) $n=6$: The 19 rational curves together with the classes
$L/3$ and $M/3$ generate a rank-19 lattice $W$ with discriminant
$2 \cdot 3^2 \cdot 5$. If $W \not= NS(Y)$, the only possibility 
would be $|NS/W|=3$. But then $|NS^{\vee}:NS| \leq 10$ and
$|NS^{\vee}/W| \leq 3 \cdot 10$. However there are the two independent
classes
$$\frac{1}{3}(N_1-N_2-N_3+N_4), \quad 
\frac{1}{3}(N_5-N_6-N_7+N_8) \quad \in NS^{\vee}$$
generating in $NS^{\vee}/W$ a subgroup isomorphic with
$\Z_3 \times \Z_3$, contradiction.  

b) $n=8$: $V$ together with the classes
$L/2$ and $L'/2$ generates a rank-19 lattice $W$ with discriminant
$2^4 \cdot 3 \cdot 7$. Here we have the three independent classes
$$M:=M_3+M_4, \quad M':=M_1+M_2+M_3, \quad
R:=R_1+2R_2+3R_3$$
with $M/2, \, M'/2, \, R/4 \in NS^{\vee}$. They generate in 
$NS^{\vee}/NS$ a subgroup isomorphic to $(\Z_2)^2 \times \Z_4$. 
This implies that $NS/W$ does not have 2-torsion. We find $NS=W$.

c) $n=12$: The lattice $V=W$ spanned by the 19 rational curves has
discriminant $2^3 \cdot 3^2 \cdot 5 \cdot 11$. If $W \not= NS$, 
then $NS/W$ would have 2-torsion or 3-torsion. However the three
classes $M_i/2, i=1,2,3, \in NS^{\vee}$ generate a subgroup 
isomorphic to $\Z_2^3$ in $NS^{\vee}/NS$. This excludes 2-torsion
in $NS/W$. And the two classes $(N_1-N_2)/3, (N_3-N_4)/3$ generate in
$NS^{\vee}/NS$ a subgroup isomorphic with $\Z_3^2$. This excludes
3-torsion in $NS/W$. We have shown: $NS(Y)=W$. \qed 

\subsection{The special cases}

Just as before we denote by $V \subset NS(X)$ the sub-lattice
spanned by the rational curves from sect. 4.1. Now it has
rank 20. In the same way, as in sect. 6.1 we check, 
that for $n=6$ the classes $L/3, L'/3$
and for $n=8$ the classes $L/2,L'/2$ in $NS(X)$ exist.
Intersecting with the twentieth rational curve we find,
that the curves can be labelled as in the diagrams of
sect. 4.2.

\begin{theo} In all cases $NS(X)$ is spanned by the classes from sect. 6.1
and the twentieth rational curve. The discriminants of the lattices are
$$ 
\begin{array}{c|cccc|cccc|cccc}
\mbox{case} & 6,1 & 6,2 & 6,3 & 6,4 & 8,1 & 8,2 & 8,3 & 8,4 &
                           12,1 & 12,2 & 12,3 & 12,4 \\ \hline
d       & -15 & -60 & -60 & -15 & -28 & -84 &-168 &-112 & 
                                -660 & -440 & -792 & -132 \\
\end{array}
$$ \end{theo}

Proof. The discriminants in the above table are those of the lattice $W$
spanned by the curves from 6.1 and by the twentieth rational curve.

{\em The case $n=6:$} The discriminant $d=-15$ is square-free. So in these two
cases $NS(X)$ does not contain a proper extension of $W$. The discriminant
$d=-60$ in cases $6,2$ and $6,3$ contains the square $2^2$. However in these
cases the class $M/2$ belongs to $NS(X)^{\vee}$, but not to $NS(X)$.
Also in these cases $NS(X)$ cannot contain a proper extension of $W$.

{\em The case $n=8$:} Posssible extensions of $W \subset NS(X)$
correspond to quadratic factors in the discriminants of the
table above. However, a possible extension in
$$ \begin{array}{c|cl}
\mbox{case} & \mbox{of degree} & \mbox{contradicts} \\ \hline
8,1 & 2 & (M_1+M_2+M_4)/2 \in NS(X)^{\vee} \\
8,2 & 2 & (M_3+M_4)/2  \in NS(X)^{\vee} \\   
8,3 & 2 & (R_1+2R_2+3R_3)/4  \in NS(X)^{\vee} \\
8,4 & 2 \mbox{ or }4 & (M_3+M_4)/2, \, (R_1+2R_2+3R_3)/4 \in NS(X)^{\vee} \\
\end{array} $$    

{\em The case $n=12:$} The possible extensions in 
$$ \begin{array}{r|cl}
\mbox{case} & \mbox{of degree} & \mbox{contradicts} \\ \hline
12,1 & 2 & M_2/2, M_3/2 \in NS(X)^{\vee} \\
12,2 & 2 & (M_1+M_2)/2, M_3/2 \in NS(X)^{\vee} \\
12,3 & 6 & M_1/2, (M_2+M_3)/2, (N_1-N_2)/3 \in NS(X)^{\vee} \\
12,4 & 2 & M_1/2, M_2/2 \in NS(X)^{\vee} \\
\end{array} $$ \qed

\section{Comments}

1) Denote by $M_k$ the moduli-space of abelian surfaces with
level-(1,k) structure 
In [Mu] the quotients $\P_3/G_6$, resp.  $\P_3/G_8$ are identified
with the Satake-compactification of $M_3$, resp $M_4$, 
and
$\P_3/G_{12}$ is shown to be birationally equivalent with the 
Satake-compactification of $M_5$. However the proof there is not
very explicit. It is desirable to have an explicit
identification of the quotient $\P_3/G_n$ with the corresponding 
moduli space. The pencil $Y_{\lambda}'$ on $\P_3/G_n$
might be useful.

2) We did not consider the quotient threefold
$\P_3/G_n$. We just identified the minimal non-singular model $Y_{\lambda}$
for each quotient $Y'_{\lambda}$.
Of course it would be desirable to have a global resolution 
of $\P_3/G_n$ 
and to view our K3-surfaces as a pencil on this smooth
threefold. One would need a particular crepant resolution of
the singularities of $\Upsilon$. Such resolutions are given
e.g. in [I, IR, Ro]. We would need a resolution, where
the behaviour of the K3-surfaces can be controlled, to identify
the partial resolutions of the four special surfaces.

3) Our quotient surfaces admit a natural involution induced
by the symmetry $C$ from [S, p. 433] normalizing $G_n$, 
but not belonging to $SL(4,\C)$. It would be interesting
to identify the quotients. 

4) By [Mo] each K3-surface with Picard number 19 admits a Nikulin-involution,
an involution with eight isolated fix-points. We do not know how to 
identify it in our cases. It cannot be the involution from 3),
because this has a curve of fix-points.
It is also not clear to us, whether this Nikulin-involution exists
globally, i.e. on the total space $\Upsilon$ of our fibration.
This Nikulin-involution is related to the existence of a sub-lattice
$E_8 \perp E_8 \subset NS(Y)$. We did not manage to identify
such a sub-lattice. 

5) It seems remarkable that the Picard group of the general surface in
a pencil of K3-surfaces can be identified so explicitly, as it is done
in sect. 6. 
It is also remarkable that the quotient K3-surfaces have Picard
number $\geq 19$. Such pencils have been studied in [Mo]
and [STZ]. We expect our surfaces to have some 
arithmetical meaning. In particular the 
prime factor $n-1=5,7,11$ in the discriminant of the Picard lattices  
draws attention. In fact, the same prime factor appears in
each polynomial $s_n, \, n=6,8,12$ from [S].
It can be found too in the cross-ratio 
$CR(\lambda_1,...,\lambda_4)$ of the four special parameters in each
pencil $X_{\lambda}$ and together with strange prime factors in the
absolute invariant $j$:
$$ \begin{array}{c|ccc}
n & 6 & 8 & 12 \\ \hline
  &   &   &     \\
CR & \D \frac{5^2}{3^2} & \D \frac{7^2}{2^4 \cdot 3} 
                        & \D \frac{11^2}{2^5 \cdot 3} \\
  &   &   &     \\
j  & \D \frac{13^3 \cdot 37^3}{2^8 \cdot 3^4 \cdot 5^4}
   & \D \frac{13^3 \cdot 181^3}{2^8 \cdot 3^2 \cdot 7^4}
   & \D \frac{12 \, 241^3}{2^{10} \cdot 3^2 \cdot 5^4 \cdot 11^4} \\
\end{array} $$

\section{References}

\noindent
\begin{tabular}{lp{15.2cm}}
B & Barth, W.: K3-surfaces with nine cusps. Geom. Dedic. 72, 171-178
               (1998) \\
BPV & Barth, W., Peters, C., Van de Ven A.: Compact Complex Surfaces.
               Ergebn. Math (3), {\bf 4}, Springer 1984 \\
C & Coxeter, H.S.M.: The product of the generators of a finite group
               generated by reflections. Duke Math. J. 18, 765-782 (1951) \\
FG & Fischer, W., Grauert, H.: Lokal triviale Familien kompakter
               komplexer Mannigfaltigkeiten. Nachr. Akad. Wiss. G\"{o}ttingen 
               II Math. Phys. Kl. (1965), 89-94 \\
I & Ito, Y.:   Crepant resolutions of trihedral singularities. 
               Proc. J. Acad. Sci. 70, Ser A, 131-136 (1994) \\
IR & Ito, Y., Reid, M.: The McKay correspondence for finite subgroups
               of $SL(3,\C)$. In: Higher Dimensional Complex Varieties.
               Proc. Conf. Trento 1994, 221-240 \\   
Mu & Mukai, S.: Moduli of abelian surfaces and regular polyhedral groups.
               In: Proc. of Moduli of Algebraic Varieties Symposium,
               Sapporo (1999) \\
Mo & Morrison, D.: On K3-surfaces with large Picard number. Inv. math. 75, 
               105-121 (1984) \\
N & Nikulin, V.: On Kummer surfaces. Izv. Akad. Nauk. SSSR, Ser. Math. 39,
               1145-1170 (1975) \\
Ra & Racah, G.: Sulla caratterizzatione delle rappresentazioni irriducibili
               dei gruppi semisemplici di Lie. Rend. Accad. Naz. Linc, 
               Sci. Fis. Mat. Nat (8), {\bf 8}, 108-112 (1950) \\  
Ro & Roan, S.S.: Minimal resolutions of gorenstein orbifolds in dimension
               three. Topology 35, 489-508 (1996) \\ 
S & Sarti, A.: Pencils of symmetric surfaces in $\P_3$. J. of Alg. 246,
               429-552 (2001) \\

\end{tabular}

\noindent
\begin{tabular}{lp{15.2cm}}
SI & Shioda, T., Inose H.: On singular K3-surfaces. In: Complex Analysis
               and Algebraic Geometry. Iwanami-Shoten, Tokyo (1977),
               117-136 \\
STZ & Sun, X., Tan, S.-L., Zuo, K.: Families of K3-surfaces over curves
               satisfying the equality of Arakelov-Yau's type and 
               modularity. Preprint (2002) \\
T & Tan, S.-L.: Cusps on some algebraic surfaces and plane curves.
               Preprint (1999) \\

\end{tabular}

\vspace*{1.0cm}

\begin{center}
{\small Mathematisches Institut, Universit\"at Erlangen, Bismarckstr. 1 1/2, 91054 Erlangen, Germany\\
e-mail:barth@mi.uni-erlangen.de\\

\vspace*{0.5cm}

FB 17 Mathematik, Universit\"at Mainz, Staudingerweg 9, 55099 Mainz, Germany\\
e-mail:sarti@mathematik.uni-mainz.de\\
}
\end{center}

\end{document}